\documentclass{amsart}
\usepackage{amsmath}
\usepackage{amsthm}
\usepackage{amssymb}
\usepackage{mathtools}
\usepackage{stmaryrd}
\usepackage{hyperref}
\usepackage{cleveref}
\usepackage{url}
\usepackage{tikz-cd}

\DeclareMathOperator{\QQQ}{\mathbb{Q}}
\DeclareMathOperator{\RRR}{\mathbb{R}}
\DeclareMathOperator{\NNN}{\mathbb{N}}
\DeclareMathOperator{\ZZZ}{\mathbb{Z}}
\DeclareMathOperator{\AAA}{\mathcal{A}}

\DeclareMathOperator{\CC}{\mathcal{C}}
\DeclareMathOperator{\DD}{\mathcal{D}}
\DeclareMathOperator{\UU}{\mathcal{U}}
\DeclareMathOperator{\VV}{\mathcal{V}}

\DeclareMathOperator{\TT}{\mathcal{T}}
\DeclareMathOperator{\PPP}{\mathcal{P}}
\DeclareMathOperator{\Obj}{Ob}
\DeclareMathOperator{\op}{^{\! \mathrm{op}}}

\DeclareMathOperator{\cl}{cl}
\DeclareMathOperator{\Mod}{Mod}
\DeclareMathOperator{\fpmod}{mod}
\DeclareMathOperator{\fp}{fp}
\DeclareMathOperator{\Fun}{Fun}
\DeclareMathOperator{\Sp}{Sp}

\DeclareMathOperator{\Spec}{Spec}
\DeclareMathOperator{\add}{add}
\DeclareMathOperator{\Add}{Add}
\DeclareMathOperator{\Mat}{Mat}
\DeclareMathOperator{\proj}{proj}

\DeclareMathOperator{\supp}{supp}
\DeclareMathOperator{\Ab}{Ab}
\DeclareMathOperator{\Img}{Im}
\DeclareMathOperator{\img}{im}
\DeclareMathOperator{\pr}{pr}

\DeclareMathOperator{\Ker}{Ker}
\DeclareMathOperator{\Coker}{Coker}
\DeclareMathOperator{\id}{id}
\DeclareMathOperator{\Tcl}{\TT^{\mathrm{cl}}}
\DeclareMathOperator{\Tid}{\TT^{{d}_I}}

\DeclareMathOperator{\idl}{Idl}
\DeclareMathOperator{\Hom}{Hom}

\newcommand{\principal}[1]{#1 \! \downarrow}
\newcommand{\frincipal}[1]{#1 \! \uparrow}
\newcommand{\qrincipal}[1]{\widehat{#1} \!   \downarrow}
\newcommand{\Ob}[1]{\Obj \left( #1 \right)}
\newcommand{\FunR}[1]{\Fun_{ #1}}

\theoremstyle{definition}
\newtheorem{defi}{Definition}
\numberwithin{defi}{section}
\theoremstyle{plain}
\newtheorem{thm}[defi]{Theorem}
\newtheorem{lem}[defi]{Lemma}
\newtheorem{cor}[defi]{Corollary}
\newtheorem{pro}[defi]{Proposition}
\theoremstyle{remark}
\newtheorem{exm}[defi]{Example}
\newtheorem{rem}[defi]{Remark}

\title{The spectrum of the real line}
\keywords{Persistence module, Ziegler spectrum, real numbers, totally ordered set, poset representation, order topology}
\subjclass{16G20 (55N31 18G05)}
\author{Jan-Paul Lerch}
\address{Jan-Paul Lerch \\ Faculty of Mathematics \\
Bielefeld University \\
PO Box 100 131\\
D-33501 Bielefeld }
\email{lerch@math.uni-bielefeld.de}
\thanks{The author has been supported by the Alexander von Humboldt Foundation in the framework of an Alexander von Humboldt Professorship endowed by the German Federal Ministry of Education and Research.}
\thanks{\emph{Date:} \today}

\begin{document}
\begin{abstract}
Motivated by the study of persistence modules over the real line, we investigate the category of linear representations of a totally ordered set.
We show that this category is locally coherent and we classify the indecomposable injective objects up to isomorphism.
These classes form the spectrum, which we show to be homeomorphic to an ordered space.  
Moreover, as the spectral category turns out to be discrete, the spectrum parametrises all injective objects. 

Finally, for the case of the real line we show that this topology refines the topology induced by the interleaving distance, which is known from persistence homology.
\end{abstract}
\maketitle
\section{Introduction}
Let $k$ be a field and let $T$ be a totally ordered set. 
The spectrum of indecomposable injective objects of a locally coherent category as a topological space was introduced in \cite{Herzog97} and \cite{HK97}, independently. 
We investigate the category $\Mod kT$ of $k$-linear representations of $T$ and show that it is locally coherent (cf. \Cref{lem:fpmodAbelian}).
The main result is the following:
\begin{thm}[{cf. \Cref{thm:main}}]
There is a homeomorphism 
\begin{align*}
\Sp \Mod kT \rightarrow \idl T
\end{align*}
between the spectrum of $\Mod kT$ and the space of ideals of $T$ in the order topology.
\end{thm} 
A locally coherent category $\AAA$ is a locally finitely presented category such that the full subcategory of finitely presented objects $\fp \AAA$ is abelian.
The spectrum $\Sp \AAA$ of $\AAA$ is then defined to be the set of isomorphism classes of indecomposable injective objects.
For the topology, see \Cref{subsec:spec}.

The closed subsets of the spectrum of indecomposable injectives of a locally coherent Grothendieck category $\AAA$ are in bijection with Serre subcategories of the full subcategory of finitely presented objects $\fp \AAA$ and with hereditary torsion pairs of finite type for $\AAA$, see \cite{HK97}.
For a Serre subcategory of $\fp \AAA$, taking the minimal subcategory closed under all direct limits yields a localising subcategory of $\AAA$.
The categories obtained in this way, which are called \emph{localising subcategories of finite type}, are consequently also completely parametrised by the closed subsets of $\Sp \AAA$.

\vspace{0.2cm}

Note that there are several other related notions of spectra.
In algebra the most prominent example of a spectrum is the the prime spectrum $\Spec A$ of a commutative ring $A$.  
Since rings themselves can get very complicated, the spectrum helps to understand the structure better. 
Particularly, its topology, the Zariski topology, is an essential tool to understand the module category of this ring.

In this context, Matlis made the observation in his thesis \cite{Mat58} that for commutative Noetherian rings the prime spectrum $\Spec A$ is in bijection with the isomorphism classes of indecomposable injective modules. 
Gabriel generalised some of Matlis' results in his thesis \cite{Ga62}, where he defined the \emph{spectrum of a Grothendieck category} $\AAA$ to be the set $\Sp \AAA$ of isomorphism classes of indecomposable injective objects, but without defining a topology, yet. 

Later, Ziegler introduced the \emph{spectrum of pure-injective modules} in his model theoretic paper \cite{Zi84}, where he defined a topology on this set. 
For a ring $R$ he defined the spectrum to be the set of isomorphism classes of indecomposable pure-injective modules, which is now known as the \emph{Ziegler spectrum}. 

For a commutative noetherian ring $A$, the spectrum $\Sp \Mod A$ with the Ziegler topology was shown to be the Hochster dual of the Zariski spectrum $\Spec A$.
See also \cite{Prest09} for an extensive treatment of the spectrum and its model theoretic and categorical aspects.
The Ziegler spectrum and the spectrum with their respective topologies were generalised to locally coherent categories in both \cite{Herzog97} and \cite{HK97}, independently. 
In fact, the Ziegler spectrum and the spectrum are related via embeddings into functor categories, see \cite{GJ81} and \cite{CB94}.

We use the representation theory of partially ordered sets and totally ordered, which was promoted in the seventies by Mitchell, as part of his more general surveys on the representation theory of small categories, see \cite{RWSO}, and by Nazarova and Ro\u{\i}ter and the Kyiv school \cite{NR72}.
It gained importance in recent years due to the omnipresence of persistence modules in persistence homology and topological data analysis.
See for example \cite{Oudot2015} for an introduction to this topic.

\vspace{0.2cm}
The aim of this paper is to determine the spectrum $\Sp T := \Sp \Mod kT$ of the category of $k$-vector space representations of a totally ordered set $T$ for any field $k$.
This is particularly motivated by the special case of representations of the real line, which generalises the classical case of representations of linearly ordered quivers of Dynkin type $A_n$.
A question of interest is how the categories of representations of a discrete quiver $A_n$ and of a continuous set $\RRR$ differ.
We will see that a discrete totally ordered set without density phenomena will have a fundamentally different spectrum than for example $\RRR$.

The first aim of this paper is to show that the spectrum $\Sp T$ is homeomorphic to the set of ideals $\idl T$, which is an order space (\Cref{thm:main}).
To achieve this, the category of finitely presented representations is determined and shown to be abelian (\Cref{sec:fp}), whereupon we observe that the spectrum is completely parametrised by the set of ideals of $T$ (\Cref{subsec:indecinj}).
We observe that the spectrum parametrises all injective objects, already (\Cref{subsec:spectralcat}).

Subsequently, the closure operator which defines the topology on the spectrum is made explicit (\Cref{pro:Topology}).
In particular, a subset of the spectrum is closed if and only if it is closed under all admissible element wise unions and intersection.  

Thereafter some properties of the topology are discussed and the topology is revealed to be equivalent to the well-known order topology on the the set of ideals (\Cref{thm:main}).

Lastly, we compare the topology in the special case $T=\RRR$ with the topology induced by the so called \emph{interleaving distance}, which is important to study the stability of persistence modules, being the prerequisite for data analysis (\Cref{sec:interl}).
It turns out that the Ziegler topology on the spectrum refines the topology induced by the interleaving distance, if one distinguished point is removed.

\subsection*{Acknowledgement}
This is part of an ongoing PhD project. I thank my advisor Henning Krause and Bill Crawley-Boevey for their support and discussions. Moreover, I am very grateful to Janina Letz for reading several versions of this paper and providing many helpful comments. I would also like to thank Ulrich Bauer and Rosanna Laking for their interest and questions on large modules and the interleaving distance (see \Cref{exm:largemod} and \Cref{sec:interl}) and on superdecomposable injectives (see \Cref{subsec:spectralcat}), respectively.

\section{Preliminaries and Notation}
\subsection{Representations of small categories} \label{subsec:smallCats}
We make extensive use of the notation in Mitchell's \cite{RWSO}, which advertises the idea that representations of small categories can basically be treated as modules over some `ring with several objects'. 
Notwithstanding, we use the terms \emph{preadditive} for $\ZZZ$-enriched categories, and the term \emph{additive} for preadditive categories which are further closed under finite coproducts.

We do representation theory with enriched categories to make use of the theory of additive functors:
for a commutative ring $k$, a \emph{$k$-category} is a $\Mod k$-enriched preadditive category $\AAA$. 
Equivalently, there is a ring homomorphism from $k$ to the natural endotransformations of the identity functor of $\AAA$.

The categories of covariant $\ZZZ$-linear functors between preadditive categories $\CC$ and $\DD$ are denoted by $\FunR \ZZZ (\CC, \DD)$.
The morphisms are the $\ZZZ$-linear natural transformations.
An \emph{additive functor} is a $\ZZZ$-linear functor which preserves finite coproducts, and the categories of such functors are denoted by $\Add(\CC, \DD)$, for $\CC$ and $\DD$ additive categories.

Let $\AAA$ be a cocomplete additive category. 
An object $X $ in $\AAA$ is \emph{finitely presented} if the functor $\Hom_{\AAA}(-, X)$ preserves filtered colimits.
For the full subcategory of finitely presented objects we write $\fp \AAA$.
The category $\AAA$ is called \emph{locally finitely presented} if $\fp \AAA$ is essentially small and every object of $\AAA$ is a filtered colimit of objects in $\fp \AAA$.

Let $\CC$ be a small $k$-category.
Then we consider the module category of $\CC$, which is $\Mod \CC : = \FunR \ZZZ (\CC, \Ab) $. 

Not every category is $k$-linear. 
But given any small category $\CC$ and any commutative ring $k$ we can take the \emph{$k$-linearisation} of $\CC$, denoted by $k\CC$, in the sense of Mitchell \cite{RWSO}. 
This is the uniquely defined $k$-linear category which has the same objects as $\CC$ and where the morphism set between two objects $x,y$ is given as the free $k$-module generated by the elements in $\Hom_{\CC}(x,z)$, such that the composition is $k$-bilinear.
Note, that the category of modules $\Mod (k \CC)$ is equivalent to $\Fun (\CC, \Mod k)$, the category of all functors from $\CC$ to the category of $k$-modules.

The category of finitely presented objects $\fpmod \CC := \fp \Mod \CC$ consists of those representations $M$ which fit into the right exact sequence
\begin{align} \label{eq:fp}
	\bigoplus_{i=1}^m \Hom_{\CC}(X_i, -) \rightarrow \bigoplus_{j=1}^n \Hom_{\CC}(Y_j,-) \rightarrow M \rightarrow 0 \ 
\end{align}  
for some $X_i, Y_j$ in $\CC$.
The representable functors $\Hom_{\CC}(X, -)$ are projective generators of $\Mod \CC$ and therefore $\Mod \CC$ is a locally finitely presented Grothendieck category.
Moreover, we write $\proj \CC$ for the full subcategory of finitely generated projective objects of $\Mod \CC$.

\begin{rem}
Note that in other places the notation $\Mod \CC$ is used for the category of additive contravariant functors $\Add (\CC \op , \Ab)$, instead of covariant functors.
The relation between both concepts is addressed after \Cref{lem:ModProj}.
\end{rem}

Now, let $M$ be a  functor in $\Mod \CC$.
We write ${ M_p:= M(p)} $ for the evaluation of $M$ at an object $p \in \CC$. 
Moreover, for $x \in M_p$ and $p \in \CC$, we just write $x \in M$ and define $|x|:=p$.
Note that $M_p$ is a $k$-module and $M(\lambda): M_p \rightarrow M_q$ is a $k$-module homomorphism for any $\lambda \colon p \rightarrow q$.

We also need the \emph{support} function of a module:
\begin{align*}
\supp \colon \Mod \CC \rightarrow \PPP \left( \Ob \CC  \right), \ \supp M = \left\lbrace c \in \CC \mid M_c \neq 0 \right\rbrace,
\end{align*}
where $\PPP(-)$ denotes the power set.

Let $h^{-} \colon \CC^{op} \rightarrow \Mod \CC$ denote the contravariant Yoneda embedding, sending $c \in \CC$ to the representable functor $h^c : = \Hom_{\CC}(c,-)$.
The category $\Mod \CC$ is a Grothendieck category and the representable functors form a generating subset of it, so there is a generating epimorphism 
\begin{align*}
\gamma:\bigoplus_{j \in J} h^{c_j} \rightarrow M,
\end{align*}
for some $c_j \in \CC, j \in J$ for an indexing set $J$.
For every $j\in J$, we fix an $x_j \in M_{c_j}$, the image of the identity element in $\Hom_{\CC}(c_j,c_j)$ under the morphism $\gamma(c_j)$. 
These elements $x_j$ associated with a generating epimorphism $\gamma$ are called \emph{generators} of $M$ and they are uniquely defined by the Yoneda Lemma.
For a subset $J' \subseteq J$ we write $\langle x_j \mid j \in J' \rangle$ for the subrepresentation generated by the images of the summands $h^{c_j}, j \in J',$ under $\gamma$.

Conversely, every element $x \in M$ defines a unique homomorphism $h^{|x|} \rightarrow \langle x \rangle$.
In this notation, an element $x \in M_c$ for $c \in \CC$ can be expressed in terms of those generators with $\langle x_j \rangle_c\neq 0$.
So, for adequate morphisms $\lambda_j \in \Hom_{\CC}({c_j}, c)$, we can write 
\begin{align*}
x= \sum_{j \in J} M(\lambda_j)( x_j) ,
\end{align*}
where only finitely many $\lambda_j$ are non-trivial. 
A subset of a set of generators is a set of generators, again, if all elements in $M$ can be expressed as a linear combination in terms of these generators, as above. 
A set of generators is called \emph{linearly independent} if there are no non-trivial linear combinations of zero in terms of them. 
Equivalently, $M$ is the direct sum of the images of $h^{c_j}$ under $\gamma$.

\subsection{The spectrum of a locally coherent category} \label{subsec:spec}
Let $\mathcal{A}$ be a locally finitely presented Grothendieck category. 
It is known that the isomorphism classes of indecomposable objects of $\AAA$ form a set, as every indecomposable injective is isomorphic to the injective envelope of a quotient of a generator of $\AAA$.
To simplify notation we define the \emph{spectrum} $\Sp \AAA$ of $\AAA$ to be any set of representatives for every isomorphism class of indecomposable injectives in $\AAA$.

If $\AAA$ is locally coherent, meaning that the category $\fp \AAA$ of finitely presented objects of $\AAA$ admits kernels (and therefore is abelian), the \emph{Ziegler topology} on $\Sp \AAA$ is given by the Kuratowski closure operator $\UU \mapsto \overline{ \left( \UU\right) } := \left( {}^{\perp} \UU) \right)^{\perp},$ where
\begin{align*}
\CC^{\perp} &= \left\lbrace X \in \Sp \AAA \mid \Hom_{\CC}(C, X)=0 \text{ for all } C \in \CC \right\rbrace \ \text{and}\\
{}^{\perp} \UU &= \left\lbrace C \in \fp \AAA \mid \Hom_{\CC}(C, X)=0 \text{ for all } X \in \UU \right\rbrace
\end{align*}
for all subsets of objects $\CC \subseteq \fp \AAA$ and $\UU \subseteq \Sp \AAA $.
See e.g.\ \cite[{Lemma 12.1.12}]{HKVol2} for this.

These orthogonalities also induce an inclusion reversing bijection
\begin{align*}
( \ ) ^ {\perp} : \left\lbrace \text{Serre subcategories of } \fp \AAA  \right\rbrace \leftrightarrow \left\lbrace \text{closed subsets of } \Sp \AAA \right \rbrace : {}^{\perp}( \ )
\end{align*}
of the lattices of Serre subcategories of $\fp \AAA$ and the local of closed subsets of the spectrum $\Sp \AAA$, which is the opposite of the category of open sets.
Note that the bijection extends to another incarnation of \emph{Stone duality} (see \cite{PJ82}), as the lattice of Serre subcategories forms a frame.
\subsection{Representations of partially ordered sets}
Let $(P, \leq )$ be any partially ordered set. 
This intrinsically carries the structure of a small category:
the objects are the elements of $P$ and there is a unique morphisms between two elements $x, y \in P$ for every relation $x \leq y$.

For the remaining paper, let $k$ be a field and let $kP$ denote the $k$-linearisation of $(P, \leq)$. 
It is clear that for $p,q$ in $P$ the $k$-dimension of $\Hom_{kP}(p,q)$ is at most 1.

Consider the category $\Mod kP$ of representations of $kP$ as in \Cref{subsec:smallCats}. 
Now it is not hard to see that a representation of $M$ is equivalent to a family of vector spaces indexed by $P$, and vector space homomorphism between them indexed by the relation $\leq$. 
In the context of persistence theory, such representations $M$ are called \emph{persistence module} and the morphisms assigned to the relations are called the \emph{structure maps of $M$}.

An \emph{interval} $I$ in $P$ is a subset which is convex and connected, meaning that for all $x,y,z \in P$ with $x \leq y \leq z$ we have $y \in I$  if $x,z \in I$, and that the Hasse diagram associated with $P$ has only one component, meaning that every two points in $P$ are connected by a not necessarily finite zig zag.
Recall that the \emph{interval representation} $k_I$ for an interval $I \subseteq P$ is the representation defined by the property that $k_I(p)=k$ exactly when $p\in I$ and $k_I(p)=0$, otherwise, and for all $p,q \in I$ such that $p \leq q$, the map $k_I(p \leq q)$ is the identity map.

If $P=T$ is a totally ordered set,
we use the familiar notation for intervals:
for $a < b \in T$ the interval $( - \infty, b ]$ (or $(- \infty, b ) )$ denotes the subset of elements which are (strictly) smaller than $b$. 
The intervals $[a, \infty)$ and $(a, \infty)$ are defined dually.
We define the intervals $(a,b), \ [a,b) \ (a,b]$ and $[a,b]$ by taking the obvious intersections.

A representation is called \emph{pointwise finite dimensional} if for all $x \in P$ we have that $\dim M_x < \infty$.
In general, there is no hope to parametrise representations of an arbitrary partially ordered set: the problem is wild, that is not classifiable even for finitely generated representations.
But for some cases, like pointwise finite representations of totally ordered sets $T$, it is possible:

\begin{thm}[{\cite[Theorem 1.2]{Botnan2018a}}] \label{thm:Barcode}
Let $k$ be a field and $T$ a totally ordered set.
Then every pointwise finite-dimensional object in $\Mod kT$ is isomorphic to a direct sum of interval representations. 
This decomposition is unique up to isomorphism. 
\end{thm}

\begin{rem}
This is the most general version of the famous \emph{barcode theorem}, so far.
If no further assumptions on the shape of the category $\CC$ are made, one still has a decomposition:
any pointwise finite dimensional representation $M$ in $\Mod \CC$ is isomorphic to a direct sum of indecomposable modules with local isomorphism ring \cite[Theorem 2.1]{Botnan2018a}, see also \cite[§ 3.6]{GR97}.

Unlike many other results cited here, the barcode theorem does not readily generalise to partially ordered sets.
The quiver $D_4$ for example already has an indecomposable representation which is not an interval module.

The requirement of pointwise finite dimension can also not be lifted, because there are big modules which are not interval decomposable. 
See the following example which is inspired by \cite[3.2.1. Beispiele]{Hoppner81}:
\end{rem}

\begin{exm} \label{exm:largemod}
Let $M = \prod_{i \in \NNN_0} k_{[-i, \infty)}$ in $\Mod k\RRR$. 
This representation is not decomposable into interval modules:
assume there is some sum of interval modules $S= \bigoplus_{j \in J} k_{I_j}^{(\nu_j)}$ with cardinals $\nu_j$ such that there is an isomorphism $\varphi \colon S \rightarrow M$. 

The structure maps of $M$ are all injective and they are the identity on the positive real numbers and on every interval $[-i, -i+1)$ for $i \in \NNN$.
Therefore, the intervals $I_j$ must all be of the shape $I_j= [j, \infty)$ for $j \in \ZZZ_{\leq 0}$ or $I_{\infty}= \RRR$.
Because $\Coker M(-i \leq -i+1) \cong k$ for all $i \in \NNN$, we get $\nu_{j} = 1$ for all $j \in \ZZZ_{\leq 0}$ by the isomorphism.
So we get 
\begin{align*}
S \cong \bigoplus_{i \in \NNN_0} k_{[-i, \infty)} \oplus k_{\RRR} ^ {( \nu_{\infty})}.
\end{align*}
It follows that $\nu_{\infty}$ must equal $2 ^{\aleph_0}$, because $\dim M_0= 2 ^{\aleph_0}$, while the left summand has only point wise countable dimension at $0$.
In particular, the right summand is nonempty.

Now choose any nonzero $x \in k_{\RRR} ^ {( \nu_{\infty})} (r)$ for some $r \in \RRR$.
Since all structure maps of this summand are identities, for every $s \in \RRR$ with $s < r$ there is a unique lift $x_s \in k_{\RRR} ^ {( \nu_{\infty})} (s)$ along the structure map: $k_{\RRR} ^ {( \nu_{\infty})}(s \leq r)(x_s)=x$.
We fix some non positive integer $s< r$.
Let $y= \varphi ^{-1} _r (x)$ and $y_s= \varphi ^{-1} _s (x_s)$.

By definition of $M$, the projection on the $p$th factor $\pr_{p} (M(t))= 0$ vanishes for all $p >t \in \ZZZ_{\leq 0} $.
So, if $p > s$ is a non positive integer, then $\pr_p(y_s)=0$.
By construction we have $ \pr_t \left( M(s \leq r) (y_s) \right) = \pr_t(y)$ for all $t \in \ZZZ_{\leq 0}$, therefore $\pr_t (y) = 0$ for all $t \in \ZZZ_{\leq 0} $ with $t > s$.
But since $s\in \ZZZ_{\leq 0}$ can be chosen arbitrarily smaller than $r$, we get $\pr_i (y)= 0 $ for all $i \in \ZZZ_{\leq 0}$ and therefore $y=0$ and thus $x=0$, leading to a contradiction.
\end{exm}

\section{Finitely presented objects in  \texorpdfstring{$\Mod kT$}{Mod kT}} \label{sec:fp}
The goal of this section is characterise the subcategory $\fpmod kT$ of finitely presented objects of $\Mod kT$ and to show that it is abelian. 

Let $k$ be a field and $P$ a partially ordered set.
It is known that all projective representations of partially ordered sets are free, meaning that they can be decomposed into a direct sum of representable functors. 
This even holds more generally for all projectives in $ \FunR \ZZZ (\ZZZ P, \AAA)$, where $\AAA$ is a (generalised) module category, see {\cite[Proposition 5]{HL81}}.

If $P=T$ is a totally ordered set, then the indecomposable projectives in $\Mod kP$ are isomorphic to the interval representations of the form $k_{[a,\infty)}$.
The structure maps of these modules are all injective.
In fact, this is an aspect of a much more general phenomenon:

\begin{thm}[{\cite[Satz 1]{Brune78}}] \label{thm:Brune}
Suppose $\AAA$ is a Grothendieck category and $P$ is a partially ordered set without cycles which is downward directed. 
Then an object $M$ in $ \Fun (P, \AAA)$ is flat if and only if all structure maps are pure-injective and $M(p)$ is flat for all $p$ in $P$. \end{thm}
Recall that in the category of (finite dimensional) vector spaces, every object is flat and the pure-injective morphisms are precisely the injective morphisms. 
So we get that a representation in $ \Mod kT$ is flat if and only if its structure maps are injective.

Using the preceding theorem and the main structure theorem for persistence modules \Cref{thm:Barcode}, we show that $\fpmod kT$ is closed under finitely generated submodules of finitely generated projectives.
Note that this is not true in general.
In fact, representations of finite quivers with admissible relations do not have this property.

\begin{lem} \label{lem:SubProj}
Every finitely generated submodule of a finitely generated projective module over $kT$ is projective. 
\begin{proof}
Let $P$ be finitely generated projective and $M \subseteq P$ a submodule having generators $x_1,...,x_n$. 
From this we can find a linearly independent generating set of $M$ in the following way:

If the given set is not linearly independent, we can express $0 \in M_p$ for some $p \in P$ as a non-trivial linear combination
\begin{align*}
0= \sum_{i = 1}^n \mu_i \cdot M(|x_i| \leq p)(x_i)
\end{align*}
with $\mu_i $ in $k$.
Among all indices with nonzero coefficient $\mu_i$, choose an index $m$, with $|x_m|$ maximal.
Then, since all structure maps are injective by \Cref{thm:Brune}, the equation above lifts uniquely to
\begin{align*}
 -\mu_m x_m = \sum_{i \neq m} \mu_i \cdot M(|x_i| \leq |x_m|)(x_i),
\end{align*}
and thus 
\begin{align*}
x_m = \frac{-1}{\mu_m} \sum_{i \neq m} \mu_i \cdot M (|x_i| \leq |x_m|)(x_i).
\end{align*}

Thus, we see that the generator $x_m$ is superfluous as it can be expressed in terms of the other generators, and therefore can be deleted.
This procedure must terminate after finitely many steps. 
Because the structure maps of $M$ are all pointwise injective the same holds for all submodules of $M$.
Thus, $\langle x_j \rangle \cong h^{|x_j|}$ for $1 \leq j \leq n$, so $M$ is a finite direct sum of representable modules. 
\end{proof}
\end{lem}

The next aim is to show that $\fpmod kT$ is abelian. 
This is necessary for the topology on the spectrum of $\Mod kT$ to be well defined.
The obstruction for this to hold is to prove the existence of kernels.
But classifying all solutions for a matrix equation in $\fpmod kT$ seems very hard, so we use an argument to reduce the problem to a simpler one:
we prove a preadditive analogue of the result, that the finitely generated projectives already determine the entire module category.
The proof uses the construction of the matrix ring category.

Let $k$ be a field and $\CC$ a small $k$-category.
The \emph{matrix ring category} of $\CC$ is the $k$-linear category $\Mat \CC$ having as objects finite sequences $(X_1,...,X_n)$ of objects of $\CC$.
A morphism to another sequence $(Y_1,...,Y_m)$ is given by a matrix having entries $m_{ij} \in \Hom_{\CC}(X_i,Y_j)$. 
There is a canonical inclusion functor ${\iota_{\CC} \colon \CC \rightarrow \Mat \CC}$, sending a morphism to a $1\times1$ matrix.

In this setting it was shown by Mitchell \cite[Lemma 1.1 et seqq.]{RWSO} that the induced functor 
\begin{align*}
\Add \left( \Mat \CC, \Ab \right) & \rightarrow \Mod \CC \\
F & \mapsto F \circ \iota_{\CC}
\end{align*}
is an equivalence of $k$-categories. 

The matrix category satisfies the universal property of the additivisation of {$k$-linear} categories.
Namely, let $\DD$ be any $k$-linear category with finite coproducts and let $F \colon \CC \rightarrow \DD$ be a $k$-linear functor. 
Then there exists a $k$-linear functor $F'\colon \Mat \CC \rightarrow \DD$ factoring $F$ as $F=F' \circ \iota_{\CC}$, which is unique up to isomorphism.

The next lemma is a preadditive analogue of a well-known result from Morita theory. 
\begin{lem} \label{lem:ModProj}
Let $\CC$ be as above.  
Then there is an equivalence ${\Mod \CC \op \cong \Mod \proj \CC}$, or with other words, $\CC \op$ and $\proj \CC$ are Morita equivalent.
\begin{proof}
Consider the the functor $F=h^{-} \colon \CC \op \rightarrow \proj \CC$, which is the contravariant Yoneda embedding. 
By the properties discussed above, the functor $h^{-} $ factors through a functor $F' \colon \Mat \left( \CC \op \right) \rightarrow \proj \CC$.
This $F'$ maps every sequence $(X_1,...,X_n)$ to the direct sum $\bigoplus_{i=1}^n h^{X_i}$ and every matrix $(m_{ij})_{ij}$ representing a morphism from $(X_1,...,X_n)$ to $(Y_1,...,Y_m) $ to a morphism represented by the matrix $(h^{m_{ij}})_{ji}$:
 
\begin{align*}
\bigoplus_{j=1}^m h^{Y_j} \xrightarrow{(h^{m_{ij}})_{ji} }  \bigoplus_{i=1}^n h^{X_i}.
\end{align*}
Applying the Yoneda Lemma on the components of the matrix reveals that the functor $F'$ is fully faithful and thus we can treat $\Mat \left( \CC \op \right) $ as a full subcategory of $\proj \CC$.

We obtain with the following sequence of functors:
\begin{align*}
\CC \op \xrightarrow{\iota_{\CC \op}} \Mat \left( \CC \op \right) \xrightarrow{F'} \proj \CC.
\end{align*}
Now, the first functor $\iota_{\CC \op}$ induces an equivalence of the module categories $\Mod \left( \CC \op \right)$ and $\Add \left( \Mat \left( \CC \op \right), \Ab \right)$ as discussed previously to this lemma. 

The functor $F'$ identifies $\Mat \left( \CC \op \right)$ with the full subcategory of $\proj \CC$ consisting of the finite direct sums of representable functors. 
The representable functors are the generators of this category, so every element in $\proj \CC$ is a retract of the image of $\Mat \left( \CC \op \right)$ considered as a subcategory of $\proj \CC$, which means this image is a \emph{cover} of $\proj \CC$. %
It then follows from \cite[Lemma 1.1 et seqq.]{RWSO} that $F'$ induces an equivalence of categories of additive functors $\Add \left( \Mat \left( \CC \op \right), \Ab \right)$ and $\Add \left( \proj \CC , \Ab \right)$.

Taken all together, we get an equivalence $\Mod \left( \CC \op \right) \cong \Add \left( \proj \CC, \Ab \right)$.
Note that the matrix category of a small additive category $\AAA$ is isomorphic to $\AAA$ itself.
Thus, $\Mat \proj \CC$ is equivalent to $ \proj \CC$, and $\Add \left( \proj \CC, \Ab \right)$ is equivalent to $\Add \left( \Mat \proj \CC, \Ab \right)$, which again is equivalent to $\Mod \proj \CC $.
\end{proof}
\end{lem}
\begin{rem}
The preceding lemma and the Matrix construction also works if $k$ is any commutative ring.
\end{rem} 

In the last step of the proof we saw that for a small additive category $\AAA$ the module category $\Mod \AAA $ is equivalent to $\Add (\AAA , \Ab)$, which is the additive definition of the module category of an additive small category.
So our notion is compatible.

Recall that a \emph{weak kernel} of a morphism $f \colon X \rightarrow Y$ is a morphism $W \rightarrow X$ which factors every morphism $g \colon Z \rightarrow X$ with $f \circ g = 0$.
So the weak kernel has the same definition of a kernel except the uniqueness.

For the remaining paper, we return to the case $\CC = kT$. 
\begin{lem} \label{lem:fpmodAbelian}
The full subcategory $\fpmod kT$ of finitely presented objects is abelian.
\begin{proof}

For an additive category $\AAA$ the full subcategory of finitely presented objects $\fpmod \AAA $ of $\Mod \AAA $ is always closed under cokernels.
Moreover, it is abelian if and only if all morphisms in $\AAA$ have weak kernels, see \cite[2.4]{CB94} and \cite[Satz 2.4]{Br70}.

If we apply  \Cref{lem:ModProj} to $\CC = \left( kT \right) \op$, we get an equivalence of $\Mod kT$ and $ \Mod \proj \left( kT \right)\op  $.
Restricting this to the subcategories of finitely presented objects, we find that it is enough to show that $\fpmod \proj \left( kT \right) \op$ is abelian.
So we only need to check that $\proj \left( kT \right) \op$ has weak kernels.

The opposite category of a $k$-linearisation of a totally ordered set is just the $k$-linearisation of the totally ordered set $T \op$, which is obtained from $T$ by reverting all relations.
For this reason, if the category $\proj kT'$ has weak kernels for every totally ordered set $T'$, then in particular $\proj \left( kT \right) \op$ has weak kernels.

In fact, $\proj kT'$ is closed under proper kernels:
for every morphism $\varphi$ of finitely generated projectives, the image $\Img \varphi$, taken in the surrounding abelian category $\Mod kT'$, is certainly finitely generated and a subrepresentation of a finitely generated projective.
By \Cref{lem:SubProj} it is projective.
Now, the image morphism $\img \varphi$ is a projection and the domain of $\varphi$ decomposes into the direct sum $\Ker \varphi \oplus \Img \varphi$. Both summands are in $\proj kT'$, as this category is closed under direct summands.
\end{proof}
\end{lem}

Next, we explicitly describe the objects in $\fpmod kT$.
\begin{lem} \label{lem:indFP}
Every finitely presented indecomposable object in $\Mod k T$ is isomorphic to an interval representation $k_{[a,b)}$ for $a,b \in T \cup \lbrace \infty \rbrace$ with $a < b$.
\begin{proof}
Let $M$ be a finitely presented indecomposable representation of $kT$. 
Then $M$ is pointwise finite dimensional and thus isomorphic to an interval module $k_J$ for some interval $J \subseteq T$, by \Cref{thm:Barcode}. 
For an indecomposable finitely generated representation of $k T$ one can find a single generator: 

Take a finite presentation 
\begin{align*}
\bigoplus_j h^{b_j} \rightarrow \bigoplus_i h^{a_i} \xrightarrow{\pi} M \rightarrow 0
\end{align*}
 for some $b_j$ and $a_i$ in $T$. 
Then the generators of the representable modules $h^{a_i}$ are mapped to some elements $x_i$ in $M$. 
But since $M$ is pointwise at most one-dimensional, they are all contained in the subrepresentation of $M$ generated by one generator with minimal support, say it has index $l$.
So, the restriction $\pi_{ \mid_{h^{a_l}}}$ is surjective. 
The kernel of $\pi$ is finitely generated and therefore projective by \Cref{lem:SubProj}.
So it must be a direct sum of representable  representations and we can assume the projective presentation to be short exact. 
Comparing the pointwise dimensions of the terms we see that $M$ must be of the form $k_{[a_l,b)}$ for $b \in T \cup \lbrace \infty \rbrace$ and $a < b$.
\end{proof}
\end{lem}

So, we see that 
\begin{align*}
	\fpmod kT = \add \left( k_{[a,b)} \mid a, \in T, b \in T \cup \lbrace \infty \rbrace \right),
\end{align*}
where $\add$ denotes the closure under finite direct sums.

\section{The indecomposable injectives} \label{sec:spec}
The purpose of this section is to determine the (indecomposable) injective objects of $\Mod kT$.

\subsection{Indecomposable injective objects} \label{subsec:indecinj}
We use the notion of ideals:
\begin{defi}
Let $T$ be a totally ordered set.
A nonempty subset $I\subseteq T$ is called \emph{ideal} if it is closed under smaller elements, meaning that whenever $x \leq y$ and $y \in I$, then also $x \in I$. 
Dually, a nonempty subset $F \subseteq T$ is called a \emph{filter} if it is closed under all greater elements.
The set of all ideals of $T$ is denoted by $\idl T$.

Let further $a$ be an element of $T$. The \emph{principal ideal generated by $a$}, denoted $\principal a$, is the smallest ideal containing $a$.
Dually, the principal filter $\frincipal a$ is the smallest filter containing $a$.
We further write $\qrincipal{a}$ for the principal ideal generated by $a$ with its generator removed: $\principal a \setminus \lbrace a \rbrace$.
\end{defi}

The set of ideals $\idl T$ is a totally ordered set, too, as it inherits this property from $T$.
Note that in this notation we have that the relation $\qrincipal x \subsetneq I$ is equivalent to $\principal x \subseteq I$, because $\principal x$ is the next greater element after $\qrincipal x$.
Similarly, $I \subsetneq \principal x$ is equivalent to $I \subseteq \qrincipal x$.

Next we classify all indecomposable injective in terms of ideals.
For this, we make use of a theorem by H\"oppner, which is stated in a slightly more general context:

\begin{thm}[{\cite[3.3.3 Theorem]{Hoppner81}}]
Let $P$ be a partially ordered set, $R$ any (not necessarily commutative) ring and consider the functor category $\Fun(P, \Mod R)$.
Given a subset $J\subseteq P$ and an $R$-module $Q$, we define $\Delta_JQ$ to be the the representation of $P$ which constantly takes the value $Q$ on $J$ and is $0$ otherwise, and which has the identity map $\id_Q$ as structure maps, whenever applicable, and the zero map, otherwise.

Then for any functor $F$ in $\Fun(P, \Mod R)$, the following are equivalent:
\begin{enumerate}
	\item $F$ is injective and indecomposable.
	\item $F \cong \Delta_JQ$ for an injective and indecomposable $R$-module $Q$ and $J\subseteq P$ a subset which is a filtered ideal. 
\end{enumerate}
\end{thm}
If we choose $P=T$ to be a totally ordered set, then the subsets of $J \subseteq T$ satisfying condition (2) are exactly the ideals of $T$.

Note that the injective indecomposables in the category $\Mod k$ are exactly the $1$-dimensional vector spaces. 
In combination with H\"oppner's Theorem, this implies that every isomorphism class of indecomposable injectives in $\Mod kT$ is uniquely determined by its support, and the support map restricts to the bijection
\begin{align*}
\supp \colon \Sp T \rightarrow \idl T.
\end{align*}
It has an inverse map $\Delta$. 
There is the canonical choice for $\Sp T$ consisting of the interval representations $k_I$ for all $I \in \idl T$, and in this case we have $\Delta (I) = k_I$. 

We can view the set $\idl T$ as a topological space with the topology being induced by the bijection $\supp$, as working with ideals instead of modules can simplify some proofs.
In this case we will also speak of the Ziegler topology on $\idl T$.

\begin{exm} \label{exm:SpR}
In the totally ordered set $\RRR$ the ideals are exactly the right closed or right open half-lines, that is the intervals $\principal{x} =(- \infty, x]$, $\qrincipal x = (- \infty, x)$ for $x \in \RRR$. 
Then the canonical choice for the spectrum is:
\begin{align*}
	\Sp \RRR = \left\lbrace k_{(- \infty, x)} \mid x \in \mathbb{R} \cup \lbrace \infty 
	\rbrace \right\rbrace \cup \left\lbrace k_{(- \infty, y]} \mid y \in \mathbb{R} \right\rbrace.
\end{align*}
The inclusion of $\idl \RRR$ into the power set $\PPP(\RRR)$ equips $\idl \RRR$ with a total order.
Furthermore, we define the double line with infinity
\begin{align*}
D_{\RRR}^{\infty}=\RRR \times \lbrace 0, 1 \rbrace \cup \lbrace (\infty, 0) \rbrace.
\end{align*}
Then we get a bijection $\idl \RRR \rightarrow D_{\RRR}^{\infty}$,
\begin{align*}
\idl \RRR & \longrightarrow D_{\RRR}^{\infty} \\
(- \infty, x) & \longmapsto (x, 0) \\
(- \infty, y] & \longmapsto (y, 1),
\end{align*}
 with $x \in \RRR \cup \lbrace \infty \rbrace$ and $y \in \RRR$.

On $D_{\RRR}^{\infty}$ we have the lexicographic order, that is the total order induced by the underlying total order of $\RRR \cup \lbrace \infty \rbrace$ and $\lbrace 0 < 1 \rbrace$. 
This means we have $(x,i) < (y,j)$ whenever $x<y$ and $(x,i)<(x,j)$ whenever $i=0$ and $j=1$.
It is clear that the map above identifying $\idl \RRR$ and $D_{\RRR}^{\infty}$ preserves the order.

For the bijection $D_{\RRR}^{\infty} \rightarrow \Sp \RRR$ compatible with the maps above we use the following notation: if we are given an element $x=(x,i) \in D_{\RRR}^{\infty}$, 
\begin{align*}
k_x := 
\begin{cases*}
k_{(-\infty, x)} \text{ if } i=0 \\
k_{(-\infty, x]} \text{ if } i=1
\end{cases*}.
\end{align*}
\end{exm}

\subsection{Morphisms to indecomposable injectives}
With the parametrisation of the spectrum at hand, we can provide a simple formula to determine if there are non trivial maps from objects of $\fpmod kT$ to elements of $\Sp T$.

\begin{lem} \label{lem:HomCondition}
Let $I $ be an ideal of $T$ and let $a<b \in T \cup \lbrace \infty \rbrace$.
Then 
\begin{align*}
\Hom_{ \Mod k T }(k_{[a,b)}, k_I)\neq 0
\end{align*}
holds if and only if $\principal a \subseteq I$ and additionally $ I \subseteq \qrincipal{b}$ if $b \neq \infty$.
In this case the dimension of this vector space is 1.

Conversely, there are no nonzero morphisms if and only if $I \subsetneq \principal a$ or $\qrincipal b \subsetneq I$ for $b \neq \infty$.
\begin{proof}
Assume that $\Hom_{ \Mod k T }(k_{[a,b)}, k_I)$ is nonzero.
Then the supports of $k_{[a,b)}$ and $k_I$ overlap: $[a,b) \cap I \neq 0$.
This is the case if and only if $\principal a \subseteq I$. 
Now, if $\qrincipal{b} \subseteq I$, then there is a point $r$ in the ideal $I$, which is not in $[a,b)$. 
In this case there are elements $q \in [a,b) \cap I$ and $r \in I \cap \frincipal b$ and there is a diagram
\begin{equation*}
\begin{tikzcd}
\left( k_{[a,b)} \right)_q = k \rar \dar["\lambda \cdot"] & 0 = \left( k_{[a,b)} \right)_r \dar \\
\left( k_I \right)_q = k \rar["\id_k"]	& k = \left( k_I \right)_r
\end{tikzcd}
\end{equation*}
for all $\lambda$ in $k$.
This diagram commutes only for $\lambda=0$, therefore there cannot be any nonzero maps. 

Conversely it is not hard to see that there is a nonzero map defined by pointwise multiplication by some unit $\lambda$ in $k$ if the condition on the ideal is satisfied.

The space of homomorphisms has dimension 1 because every homomorphism in this context is completely defined by its restriction to a point, if it is nonzero.
\end{proof}
\end{lem}

\subsection{Injective objects and the spectral category} \label{subsec:spectralcat}
As mentioned in the introduction, the Zariski spectrum parametrises all indecomposable injective modules of a commutative noetherian ring. 
Gabriel and Oberst generalised this observation with their construction of the \emph{spectral category} in \cite{GO66} (see also \cite[§6]{St75}). 
This category is obtained by localising a Grothendieck category $\AAA$ at its essential monomorphisms, which identifies every object with its injective hull.
In consequence, the isomorphism classes of objects in the spectral category are in one to one correspondence with the isomorphism classes of injective objects in $\AAA$. 
Moreover, every object in the spectral category decomposes into the sum of a semisimple object, the \emph{discrete} part, and a summand without simple subobjects, the \emph{continuous} part.
Thus, every injective object in $\AAA$ decomposes into the injective hull of a direct sum of indecomposable injectives and a \emph{superdecomposable} summand without indecomposable summands.

The following theorem by H\"oppner implies that there are no non-trivial superdecomposable injective objects in $\Mod k\RRR$. 
Recall that a representation is called \emph{uniform} if every non-trivial subrepresentation is essential, i.e.~ intersects with every other subrepresentation non-trivially, and that the injective hull of a uniform representation is indecomposable.

\begin{thm}[{\cite[3.3.6 Satz]{Hoppner81}}]
Let $R$ be any ring and $P$ a partially ordered set. 
Then every representation in $\FunR \ZZZ (P, \Mod R )$ has a uniform subobject if and only if the following hold:
\begin{enumerate}
\item Every non-trivial $R$-module has a uniform submodule.
\item The partially ordered set $P$ does not contain the binary tree which is inductively obtained by starting with a single point and adding two distinct larger elements for every point.
\item Every non-trivial representation has a non-trivial subrepresentation of the form $\Delta_J M / \Delta_L M$, where $M$ is an $R$-module, $J = \lbrace p \in P \mid p \geq i \rbrace $ for some $i \in P$ and $L \subseteq J$ a subset closed under greater elements.  
\end{enumerate}
\end{thm}
The category of vector spaces over a field is semisimple, a linearly ordered set does not contain any branches and the non-trivial image of a generating indecomposable representable object satisfies the third condition. So:
\begin{cor}
Every injective object in $\Mod kT$ is isomorphic to the injective envelope of a direct sum of indecomposable injectives. 
Equivalently, the spectral category of $\Mod kT$ is discrete.
\end{cor}
This means that all, even the large, injectives can be parametrised in terms of the spectrum $\Sp T$.

\section{The ideals of a totally ordered set}

As we have already seen in \Cref{exm:SpR}, the set of ideals $\idl \RRR$ of $\RRR$ inherits the structure of a partially ordered set from its embedding into $\PPP (\RRR)$, the power set of $\RRR$.
This holds more generally for every totally ordered set $T$.
In fact, the set of ideals is totally ordered by inclusion, as every two non-empty ideals intersect non trivially and then one already contains the other. 

There is a canonical order preserving embedding $T \rightarrow \idl T$, sending each element $a$ of $T$ to the principal ideal $\principal a$.
Moreover, the set of ideals $\idl T$ can be viewed as a completion of a totally ordered set in the sense that every subset of $\idl T$ has a supremum, even if this does not hold for $T$ itself.

\begin{exm}
For $T = \NNN$, the ideals are all principal ideals $\principal n$ for $n \in \NNN$ and $\NNN$ itself.
Therefore $\idl \NNN \xrightarrow{\sim} \NNN \cup \lbrace \infty \rbrace$ as totally ordered sets.
Note that the maximal ideal $\NNN$ is neither of type $\principal n$ or $\qrincipal n$ for any $n \in \NNN$.

For $T= \RRR$, the supremum of the set of real ideals $(-\infty, -\frac{1}{n}]$ for all natural numbers $n$ is $\qrincipal 0 =(-\infty, 0)$, which is not a principal ideal but can still be parametrised in terms of $T$.

For  $T= \QQQ$, there are even bounded subsets of $\idl T$ for which the supremum does not exist in $T$: 
let $S$ be the subset of rational numbers $x\in \QQQ$, such that $x^2 <2 $.
Then $S$ has no supremum in $T$, but the supremum of this set is $\sqrt{2}$ if embedded into $\RRR$.
However, embedded into the set of ideals $\idl \QQQ$, $S$ has a supremum, which equals the interval $(-\infty, \sqrt{2}) \cap \QQQ$.

This `completion' is not idempotent, though, and, as insinuated in \Cref{exm:SpR}, it is much more than the completion of $\QQQ$ by Dedekind-Cuts or in the Euclidean metric, which is $\RRR$.
Namely, $\idl \QQQ$ is parametrised by a copy of $\RRR$ for the bounded above non-principal ideals, a copy of $\QQQ$ for the principal ideals and a point at infinity for the entire set. 
\end{exm}

In the preceding example there are three different types of ideals: principal ideals, those which are strict lower sets and those which are neither. 
More generally, every ideal in any totally ordered set $T$ is of exactly one of the following types:
\begin{defi}
An ideal $I \in \idl T$ is of
\begin{itemize}
\item[--] Type 1, if I is a principal ideal: $I= \principal x, x \in T$.
\item[--] Type 2, if I is not a principal ideal and has a supremum in $T$.
\item[--] Type 3, if I is not a principal ideal and has \emph{no} supremum in $T$.
\end{itemize}
\end{defi}
It is immediate that an ideal of type 2 but not of type 3 is of the form $\qrincipal{a}$.
In the examples above, the ideal $(-\infty, 0)$ is of type 2 and $\RRR$ is of type 3 in $\idl \RRR$, while $\left(-\infty, \sqrt 2 \right) \cap \QQQ$ is of type 3 in $\idl \QQQ$. 

\begin{rem}
In the lattice theoretic sense, the compact objects in $\idl T$ are exactly the principal ideals $\principal a $ for $a \in T$.
So ideals of type 2 and 3 are the non-compact case.
\end{rem}

In the following, let $A ^{\complement} = X \setminus A$ denote the set complement for a subset $A \subseteq X$.
\begin{lem} \label{lem:types}
Let $I \in \idl T$ be not a principal ideal. Then:
\begin{enumerate}
\item $I$ is of type 2 if and only if $I \subsetneq \bigcap_{x \in I^{\complement}} \principal x$.
\item $I$ is of type 3 if and only if $I = \bigcap_{x \in I^{\complement}} \principal x$.
\end{enumerate} 
\begin{proof}
First, we note that for every $I\in \idl T$ we have that
\begin{align*}
I \subseteq \bigcap_{x \in I ^{\complement}} \principal x \ .
\end{align*}
Suppose that $I$ is of type 2. 
Then $I^{\complement}$ has as minimal element the supremum of $I$, which is $a$ for $I = \qrincipal{a}$. 
Thus, the intersection yields $\principal a$ and the inclusion is proper. 

Suppose $I$ is of type 3 and assume that there is some $y$ in the intersection which is not in $I$. 
Then it is not minimal with this property, since otherwise $y$ would be the supremum of $I$ and therefore $I$ would be of type 2.
So, there exists an element $z\in I^{\complement}$ with $z < y$. 
Thus, $y$ is not contained in $\principal z$ and therefore cannot be in the intersection. 

Conversely, if equality holds, then $I$ cannot have a supremum in $T$ because otherwise it would be principal.
So, it is of type 3.
\end{proof}

\end{lem}

\section{Characterisation of the closed sets}
In this section the closure operator on the spectrum $\Sp T$ is characterised. 
First, several examples for Ziegler closed sets are discussed. 
Based on these and some general properties of closures, the closed sets are then characterised.

As mentioned in \Cref{sec:spec}, it can be easier to work with ideals than with modules, so we use the identification via the support map here.
\subsection{Closed intervals in \texorpdfstring{$\idl T$}{Idl T}}
Since $\idl T$ is a totally ordered set, we can consider intervals of $\idl T$ and start with calculating the closures for some standard intervals.

\begin{lem} \label{lem:closeInt} 
Any point in $\idl T$ is (Ziegler) closed. 
For all $x \in T$ the intervals $\left[ \principal x, \infty \right)$ and $\left(-\infty, \qrincipal{x}\right]$ are closed.  
\begin{proof} 
The interval $\left[ \principal x, \infty \right)$ is closed: 
its left orthogonal consists of these sums of finitely presented representations $k_{[a,b)}$ from which there are no nonzero maps to injectives $k_I$ with $I$ in $\left[ \principal x, \infty \right)$.
By \Cref{lem:HomCondition} this holds for some $a < b$ in $T \cup \lbrace \infty \rbrace$ if and only if 
\begin{enumerate}
\item $I \subsetneq \principal a$, or
\item $\qrincipal b \subsetneq I$, if $b \neq \infty$,
\end{enumerate}
for all ideals $I$ in $\left[ \principal x, \infty \right)$.
Suppose there are such $a < b$ satisfying the first condition, namely that $I \subsetneq \principal a$ for all $I$ such that $\principal x \subseteq I$.
Then, in particular, we get $\principal x \subsetneq \principal a$, or equivalently $x < a$.
We also get $\principal y \subsetneq \principal a$ for all $y \geq x$, so in particular for $y=a$, and therefore $\principal a \subsetneq \principal a$, which is a contradiction.

The second case is equivalent to the condition $\qrincipal b \subsetneq \principal x$. But this is the case if and only if $b \leq x$. 
So:
\begin{align*}
{}^{\perp}\Delta \left( \left[ \principal x, \infty \right) \right)= \add \left( k_{[a,b)} \mid a,b \in T: a < b \leq x \right).
\end{align*}

Next we calculate the right orthogonal of this. 
It consists of all $k_I$ such that $I \subsetneq \principal a$ or $\qrincipal b \subsetneq I$ for all $a < b \leq x$.
There are no ideals $I$ satisfying the first condition, as no ideal is properly contained in $\principal m$, if $m$ is an arbitrarily small element of $T$ or a minimum, if it exists.

The second case holds if and only if $\qrincipal x \subsetneq I$, which holds exactly if $\principal x \subseteq I$ 
So we get that the closure of $\left[ \principal x, \infty \right)$ is itself.

Next, we show that $\left( - \infty, \qrincipal x \right]$ is closed:
for the left orthogonal we get all these pairs $a,b$ with $a< b \in T \cup \lbrace \infty \rbrace$ such that $I \subsetneq \principal a$ or, if $b \neq \infty$, $\qrincipal b \subsetneq I$ for all $I \subseteq \qrincipal x$.
The first condition is equivalent to $\qrincipal x \subsetneq \principal a$, wich is again equivalent to $x \leq a$.
The second condition $\qrincipal b \subsetneq I$ for all $I \subseteq \qrincipal x$ implies $\qrincipal b \subsetneq \qrincipal x$ and thus $\qrincipal b \subseteq \qrincipal x \subsetneq \qrincipal b$, hence such $a$ and $b$ do not exist.
Thus
\begin{align*}
{}^{\perp}\Delta \left( \left(-\infty, \qrincipal{x} \right] \right)= \add \left( k_{[a,b)} \mid a,b \in T \cup \lbrace \infty \rbrace : x \leq a < b \right).
\end{align*}
For the right orthogonal of this we obtain the indecomposable injectives parametrised by the ideals $I \in \idl T$ such that for all $a$ and $ b$ as above either $I \subsetneq \principal a$ or $\qrincipal b \subsetneq I$ if $b \neq \infty$.
The first condition is equivalent to $I \subsetneq \principal x$, which again is equivalent to $I \subseteq \qrincipal x$, 
The second condition is never satisfied, because $b = \infty$ is allowed here.
So the closure consists of exactly those ideals satisfying $I \subseteq \qrincipal x$ and thus $\left(-\infty, \qrincipal{x} \right]$ is closed. 

To see that a point $k_I$ is closed in $\Sp T$, we use the same technique and obtain
\begin{align*}
{}^{\perp}\lbrace k_I \rbrace = \add \left( \lbrace k_{[a,b)} \mid a < b , b \in I \rbrace \cup \lbrace k_{[a,b)} \mid a < b \in T \cup \lbrace \infty \rbrace, a \in I^{ \complement} \rbrace \right).
\end{align*}
Then 
\begin{align*}
\overline{\lbrace k_I \rbrace }
 = \add \left( \lbrace k_{[a,b)} \mid a < b , b \in I \rbrace \right)  ^ {\perp} \cap \add \left( \lbrace k_{[a,b)} \mid a < b \in T \cup \lbrace \infty \rbrace, a \in I^{ \complement} \rbrace \right) ^ {\perp}.
\end{align*}
The left term of the intersection is parametrised by all ideals $J$ with $J \subsetneq \principal a$ or $\qrincipal b \subsetneq J$ for all $a < b \in I$.
Thus, the first property vanishes and the second is equivalent to $I \subseteq J$.

The right term in the intersection is parametrised by all ideals $J$ with $J \subsetneq \principal a$ or $\qrincipal b \subsetneq J$ for all $a< b \in T \cup \lbrace \infty \rbrace$ and $a \in I^{ \complement}$. 
The second condition vanishes because $b$ can take the value $\infty$, while the first condition is equivalent to $J \subseteq I$.
It follows that $\overline { \lbrace I \rbrace } = \lbrace I \rbrace$.
\end{proof}
\end{lem}

\begin{cor} \label{cor:closedInt}
For $\idl T$ in the Ziegler topology we have the following examples of closed and corresponding open sets for all $x\in T$ and $I, J \in \idl T$, where $I$ is any ideal of type~3:
\begingroup
\renewcommand*{\arraystretch}{1.5}
\begin{align*}
\begin{array}{l | c | c}
& \mathrm{Closed \ interval} & \mathrm{Open \ complement} \\ \hline
1. & \left(-\infty, \qrincipal{x} \right]  & \left[ \principal x, \infty \right) \\
2. & \left[ \principal x, \infty \right)  & \left( - \infty, \qrincipal{x} \right] \\
3. & \left\lbrace J \right\rbrace & \left\lbrace J \right\rbrace^{\complement} \\
4. & \left(-\infty, \principal x \right]  & \left( \principal x, \infty \right) \\
5. & \left[\qrincipal{x}, \infty \right)  & \left( - \infty, \qrincipal{x} \right) \\
6. & \left(-\infty, I \right]  & \left( I, \infty \right) \\
7. & \left[I, \infty \right)  & \left( - \infty, I \right) \\
\end{array}
\end{align*}
\endgroup
\begin{proof}
\begin{enumerate}
\item[1.-3.] were covered in the previous Lemma.
\item[4.] We have $\left( -\infty , \principal x \right]  = \left(-\infty, \qrincipal{x} \right] \cup \left\lbrace \principal x \right\rbrace$, which is closed by the cases 1 and~3.
\item[5.] We have $\left[\qrincipal{x}, \infty \right) = \left\lbrace \qrincipal{x} \right\rbrace \cup \left[\principal x, \infty \right)$, which is closed by the cases 2 and 3.
\item[6.] We can write $\left(-\infty, I \right]= \bigcap_{x \in I} \left( -\infty,  \principal x \right]$, which is an intersection of closed intervals, see case 1, and therefore closed. 
\item[7.] We can write $\left[ I, \infty \right)= \bigcap_{x \in I} \left[ \principal x, \infty \right)$, which is an intersection of closed intervals, see case 2, and therefore closed. 
\end{enumerate}
\end{proof}
\end{cor}

\subsection{Topological Closures}
The following topological lemmata will be used for the characterisation of the closure.
\begin{lem} \label{lem:topological}
	Let $X$ be a topological space, $\mathcal O$ an open cover of $X$ and $S \subseteq X$. Then $S$ is closed in $X$ if and only if $S \cap U$ is closed in $U$ for all $U \in \mathcal{O}$.
\end{lem}

We use following notation for the closure of some set $S$ in $X$
\begin{align*}
	\cl_X(S)= \bigcap_{\substack{C \text{ closed} \\ S \subseteq C}} C \ .
\end{align*}

\begin{lem} \label{lem:Closure}
	Let $X$ be a topological space with open cover $\mathcal{O}$ and $S \subseteq X$ any subset. Then
	\begin{align*}
		\cl_X(S)= \bigcup_{U \in \mathcal O} \cl_U(S \cap U) \ .
	\end{align*}
\end{lem}
\subsection{Characterisation of the closure of \texorpdfstring{$\Sp T$}{Sp T}}
For the following proposition note that for a non empty subset of elements of $\idl T$ their union is in $\idl T$ and that their intersection is in $\idl T$ if and only if this subset has a lower bound in $T$.
\begin{pro} \label{pro:Topology}
Let $C \subseteq \idl T$. 
Then $C$ is (Ziegler) closed in $\idl T$ if and only if any ideal of the form $I = \bigcap_{J \in V} J$ or $ I = \bigcup_{J \in V} J$ for any non empty $V \subseteq C$ is contained in $C$.  
\end{pro}

The idea of the proof is first to show that a closed set is closed under these admissible intersections and unions, and then to show that there are no other points in the closure.
For the second step the idea is to differentiate between the types of ideals and use that ideals of type 1 (type 2) are not `reached' by unions (intersections, respectively), while ideals of type 3 can be both, an intersection and a union.
For an ideal which is neither an intersection nor a union of elements of $C$, we then choose open intervals which are part of a disjoint open cover of $\idl T$ and do not intersect $C$. 
Applying the closure operator locally we find that such ideals cannot lie in the closure.
For convenience we use the notation $\qrincipal {\infty} = T$.

\begin{proof} \label{prf:central}
Let $C = \overline C$.
We begin the proof with a reformulation of the condition
\begin{align*}
I \in C 
	& \Leftrightarrow  \Hom_{\Mod kT} (k_{[a,b)}, k_I)=0 \text{ for all } k_{[a,b)} \in {}^{\perp} \Delta (C) \\
	& \Leftrightarrow  \left( \Hom_{\Mod kT} (k_{[a,b)}, k_I) \neq 0 \Rightarrow k_{[a,b)} \notin {}^{\perp} \Delta (C) \right) \\
	& \Leftrightarrow  \left( \principal a \subseteq I \subseteq \qrincipal b  \Rightarrow k_{[a,b)} \notin {}^{\perp} \Delta (C) \right),
\end{align*}
where the last equivalence comes from \Cref{lem:HomCondition}. 

In the first step we show that a closed subset is closed under these intersections and unions and thereafter we show that no other points are added by the closure operator.

Now, let $I=\bigcap_{J \in V}J$ for some subset $V \subseteq C$. 
For any $a\in T$ with $\principal a \subseteq I$ we have that $\principal a \subseteq J$ for all $J \in V \subseteq C$.
Thus, $k_{[ a, \infty)}$ cannot be an element of ${}^{\perp} \Delta (C)$.
Similarly, let $a, b \in T$ satisfy $\principal a \subseteq \bigcap_{J \in V}J \subseteq \qrincipal{b}$.
Then there is some $J \in V$ with $J \subseteq \qrincipal b$, so we have $\principal a \subseteq J \subseteq \qrincipal{b}$.
But this implies that $k_{[a,b)}$ cannot be in ${}^{\perp} \Delta (C)$.
Thus $I$ is in $C$.

Let now $I=\bigcup_{J \in V} J$ for some $V \subseteq C$. 
Then for any $a \in T$ with $\principal a \subseteq I$, there is some $J\in V$ with $\principal a \subseteq J$.
Thus, $k_{[ a, \infty)}$ cannot be in ${}^{\perp} \Delta (C)$.
Moreover, let $a, b \in T$ satisfy $\principal a \subseteq I \subseteq \qrincipal{b}$. 
Then there is some $J \in V \subseteq C$ with $\principal a \subseteq J$.
Therefore, we get $\principal a \subseteq J \subseteq \qrincipal{b}$, hence $k_{[a,b)}$ cannot be in ${}^{\perp} \Delta (C)$.
Thus $I$ is in $C$.

It remains to show the converse: 
for all $I \in \idl T$ such that $I$ is not a union or non empty intersection of ideals in $C$, we have that $I$ is not contained in the closure $\overline{C}$.
Since every element in $C$ trivially is a union or intersection and $C \subseteq \overline{C}$, we may assume that $I$ is not in $C$.
For this, we differentiate between the three types of ideals $I$.

\textit{If $I$ is of type 1}, so $I= \principal a$ for some $a$ in $T$, then it contains its supremum and therefore cannot be a union of strictly smaller ideals.
If $a$ is not maximal in $T$, there is some $x \in I^{\complement}$, such that there is no ideal $J \in C$ with the property that $I \subsetneq J \subsetneq \principal x$:
let $V := \lbrace J \in C \mid I \subsetneq J \rbrace$, then $I \subsetneq \bigcap_{J \in V}J$ by assumption.
So we can take some $x \in \left( \bigcap_{J \in V}J \right) \setminus I$.

Now consider the interval $\left[ \principal a, \qrincipal{x} \right]$, which is an intersection of two sets that are both open and closed by \Cref{cor:closedInt}.
For $a$ maximal in $T$, choose the interval $\lbrace \principal a \rbrace$ instead.
The interval certainly contains the ideal $I$, but does not intersect $C$.
This interval together with its complement is an open covering of $\idl T$ and therefore closures can be calculated locally in this interval by \Cref{lem:Closure}.
But the closure of the empty set is the empty set, therefore $I$ is not contained in $C$.

\textit{If $I$ is of type 2}, so $I = \qrincipal a$ for some $a \in T$, then it cannot be a proper intersection of elements in $C$.
Again, we find some $y \in I$ such that there is no $J \in C$ satisfying the property $\principal y \subsetneq J \subsetneq I$:
let $V := \lbrace J \in C \mid J \subsetneq I \rbrace$, then $ \bigcup_{J \in V}J \subsetneq I$ by assumption. 
So we can take $y \in I \setminus \left( \bigcap_{J \in V}J \right)$.

Now, we can take the interval $\left[\principal y, I \right]$, which is closed, open and by construction has no intersection with $C$, and proceed as for type 1.

\textit{If $I$ is of type 3}, then it can be a union or an intersection or both by \Cref{lem:types}. 
In the first case, we can proceed like for type 1, in the second case like for type 2 and in the third case we can combine both:
we find elements $x, y \in T$ and get a closed and open interval $\left[ \principal x, \qrincipal{y} \right]$, which does not intersect $C$. 
Again, the proof is completed as for type 1. 
\end{proof}

This translates immediately to
\begin{cor}
Let $\UU \subseteq \Sp T$. 
Then $\UU$ is closed in $\Sp T$ if and only if for any non empty $\VV \subseteq \UU$ the intersection $\bigcap_{M \in \VV} \supp M$ and the union $\bigcup_{M \in \VV} \supp M$ are in $\supp \UU$ when nonempty. \qed
\end{cor}

\section{Comparison with the order topology}
\subsection{The order topology}
As discussed in \Cref{sec:spec}, there is a total order on $\idl T$. 
Since we have the support map from $\Sp T$ to $\idl T$ and because a total order induces the order topology, it is natural to ask the question, whether the support map is continuous for this topology, or how both topologies are related.
In fact, it is a homeomorphism as shown in the following.

\begin{defi}
Let $T$ be a totally ordered set. 
Then the \emph{order topology on $T$} is the topology generated by the subbasis given by the \emph{open rays} $\lbrace b \in T \mid b> a \rbrace$ and $\lbrace b \in T \mid b< a \rbrace$ for all $a \in T$.
\end{defi}

\begin{thm} \label{thm:main}
The support map $\supp \colon \Sp T \rightarrow \idl T$ is a homeomorphism, when $\idl T$ is considered to have the order topology. 
\begin{proof}
The subbasis of the order topology on $\idl T$ consists of the open rays.
The image of these open sets under $\Delta$ is open in $\Sp T$, see \Cref{cor:closedInt}.

Conversely, let $U \subseteq \idl T$ be any Ziegler open set, so $U = C ^{\complement} $ for $C : =  \supp \overline \UU $ and some $\UU \subseteq \Sp T$.
We claim that every $I \in U$ has a neighbourhood $V \subseteq U$ open in the order topology.
Suppose the contrary holds and that there is some $I \in U$ which has no neighbourhood $V \subseteq U$ open in the order topology.
This is the case if and only if for all such neighbourhoods $V$ of $I$ we have $V \cap C \neq \emptyset$.
As in the proof of \Cref{pro:Topology} we proceed with a case distinction about the types of ideals.

If $I=\principal a$ is of type 1, then every non-empty interval $\left[ \principal a, \infty \right)$ or $\left[\principal a, \principal b \right)$ for some $b>a$ contains some $J \in C$. 
This yields a set of ideals $W \subseteq C$ with $\bigcap_{J \in W}J=I$, so $I \in C$ by \Cref{pro:Topology}.
Thus we get a contradiction.

If $I= \qrincipal{a}$ is of type 2, then $I$ is union of smaller ideals in $C$: 
by assumption, every non-empty open interval $V = \left( \principal b, \qrincipal{a} \right]$ for some $b \in T$ contains some $J \in C$.
But this implies:
\begin{align*}
\bigcup_{J \in C \cap V} J = \bigcup_{\substack{J' \in \idl T \\ J' \subsetneq I}}J' = I.
\end{align*}
As before, this yields $I \in C$, which is a contradiction! 

Lastly, if $I$ is of type 3 and $I \neq T$, we can combine both previous cases and observe that for every $a, b \in T$ with $\principal a \subsetneq I \subsetneq \principal b$, the interval $\left( \principal a, \principal b \right) $ contains an element $J \in C$. 
If $I=T$ we can take an interval $\left( \principal a, T \right]$, instead.
This yields either a family of smaller ideals of which $I$ is the union of, or a family of greater elements of which $I$ is the intersection, or both.
Any way, we get $I \in C$ and the contradiction, again.

This implies that the support map induces a bijection of open sets and finishes the proof.
\end{proof}
\end{thm}

\subsection{Properties of the topology}
From the equivalence to an order topology and from the study of the closure we can obtain several properties of the spectrum.  
\begin{cor}
The space $\Sp T$ is Hausdorff.
\begin{proof}
This follows from \Cref{lem:closeInt} and \Cref{cor:closedInt} or from the equivalence in \Cref{thm:main}.
\end{proof}
\end{cor}
For the special case $T = \RRR$ we can furthermore say:
\begin{cor}
The Ziegler topology on $\Sp \RRR$ is not discrete, that is not every one point is open. \qed
\end{cor}
This phenomenon is not visible in case of for example finite linearly ordered sets, which are the same as linear$A_n$ quivers, because then the spectrum is the discrete topological space.

\begin{cor}
The space $\Sp \RRR$ is not compact in the Ziegler topology.
\begin{proof}
We have an open cover (which is even closed):
\begin{align*}
\Sp \RRR = \bigcup_{n \leq 0} \Delta \left( \left[ \principal {( n-1)}, \qrincipal{n} \right] \right) \ \cup \Delta \left( \left[ \principal 0 , \RRR \right] \right).
\end{align*}
This is a disjoint union of infinitely many sets and there certainly is no finite sub cover.
\end{proof}
\end{cor}
\begin{rem}
In the context of duality as mentioned in the introduction one is tempted to ask if one can find an interesting duality here, as well.
But taking the Hochster dual requires the underlying space to be spectral.
Spectral spaces are compact sober spaces, that is that every irreducible closed set is the closure of a unique point.
The latter property holds for $\Sp \RRR$ as it is Hausdorff.
But it is not compact, so the known approach fails to work in this context.
\end{rem}

\section{Interleaving distance} \label{sec:interl}
There is another topology known on the category of representations of the real line: the interleaving distance, see for example \cite[Definition 3.3]{Oudot2015}.
For any real number $\varepsilon$, let $[\varepsilon]$ denote the \emph{shift functor}
\begin{align*}
\Mod k \RRR & \rightarrow \Mod k \RRR \\
M &\mapsto M[\varepsilon] \\
f & \mapsto f[\varepsilon],
\end{align*}
where $M[\varepsilon](r)=M(r + \varepsilon)$ and $f[\varepsilon ] \colon M[\varepsilon] \rightarrow N [\varepsilon], f[\varepsilon](r \leq s)=f(r + \varepsilon \leq s + \varepsilon)$ for all morphisms of representations $f \colon M \rightarrow N$. 

Moreover, for any $k\RRR$-module $M$ and any $\varepsilon \geq 0$, let $\id_M^{\varepsilon}\colon M \rightarrow M [\varepsilon]$ denote the morphism of representations of $\RRR$ defined by ${\id_M^{\varepsilon}(x)=M(|x| \leq |x| + \varepsilon)(x)}$, for $x \in M$. 
\begin{defi}
An \emph{$\varepsilon$-interleaving} between two $k\RRR$-modules $M$ and $N$ is a pair of maps $\phi \colon M \rightarrow N[\varepsilon]$ and $\psi \colon N \rightarrow M[\varepsilon]$, such that $\psi [\varepsilon] \circ \phi = \id^{2 \varepsilon}_M$ and $\phi [\varepsilon] \circ \psi = \id^{2 \varepsilon}_N$.
We call $M$ and $N$ \emph{$\varepsilon$-interleaved} if there is an interleaving between them. 
The \emph{interleaving distance} of $M$ and $N$ is 
\begin{align*}
d_I(M,N)= \inf \left\lbrace \varepsilon \geq 0 \mid M \text{ and } N \text{ are } \varepsilon \text{-interleaved} \right\rbrace.
\end{align*}
\end{defi}
The definition implies that the interleaving distance assumes the value $\infty$ if there is no $\varepsilon$-interleaving for any nonnegative $\varepsilon$. 
This is one obstruction for being a metric. 
In fact, $d_i$ is an extended pseudometric; see for example the discussion after \cite[Definition 3.3 (rephrased)]{Oudot2015}.

It immediately follows that there is another topology on $\Sp \RRR$ induced by $d_I$. 
It is also clear that it is not Hausdorff (T2) if $d_I$ is a pseudometric when restricted to $\Sp \RRR$.
To find a basis for this topological space, we calculate the interleaving distance on $\Sp \RRR$ explicitly.

\begin{lem}
Let $M=k_I$ and $N=k_J$ for $I, J \in \Sp T$ with $I=\principal x$ or $I = \qrincipal x$ and $J= \principal y$ or $J= \qrincipal y$. 
Then $d_I(M, N)= |x-y|$.
\begin{proof}
We derive a criterion for the existence of maps $\phi \colon M \rightarrow N[\varepsilon]$ and $\psi \colon N \rightarrow M[\varepsilon]$, such that $\psi [\varepsilon] \circ \phi = \id^{2 \varepsilon}_M$ and $\phi [\varepsilon] \circ \psi = \id^{2 \varepsilon}_N$.
Since the defining intervals of the interval modules $M$ and $N$ are left unbounded, the maps $\id^{\delta}_M$ and $\id^{\delta}_N$ are non-trivial for every $\delta\geq 0$.
This implies that the maps $\phi, \phi[\varepsilon], \psi$ and $\psi[\varepsilon]$ are non-trivial.
Note that $\phi \neq 0 $ if and only if $\phi [\varepsilon] \neq 0$, and similarly for $\psi$. 

Since the structure maps are all $0$ or the identity, every map between $M$, $N$ or their shifts can be expressed as a pointwise multiplication by some element of $k$.
So, if there is a nonzero map between them, there also is a (nonzero) map defined by pointwise multiplication by $1$.
Then, given nonzero maps $\phi \colon M \rightarrow N[\varepsilon]$ and $\psi \colon N \rightarrow M[\varepsilon]$ which are pointwise multiplication by $1$, the concatenations $\psi [\varepsilon] \circ \phi \colon M \rightarrow M [2 \varepsilon]$ and $\phi [\varepsilon]\circ \psi \colon N \rightarrow N [2\varepsilon] $ are represented by pointwise multiplication by $1$ and therefore coincide with the maps $\id_M^{2\varepsilon}$, respectively $\id_N^{2\varepsilon}$.

So it is enough to determine the parameters $\varepsilon$ for which such nonzero maps $\varphi$ and $\psi$ exist. 
Given the case that $x < y$, then on one hand there is a nonzero map $M \rightarrow N[\varepsilon]$ if $\varepsilon >y-x$, but not if $\varepsilon < y-x$, by \Cref{lem:HomCondition}.
On the other hand there is always a map $N \rightarrow M [\varepsilon]$ for any $\varepsilon \geq 0$.
Thus, we get an $\varepsilon$-interleaving whenever $\varepsilon > y-x$, but not if $\varepsilon < y-x $.
So the interleaving distance of $M$ and $N$ is:
\begin{align*}
d_I(M,N)= \inf \left\lbrace \varepsilon\geq 0 \mid \varepsilon > y-x \right\rbrace = y-x.
\end{align*}
Note that we did not cover the case $\varepsilon = y-x$ because it would not change the infimum, but it would require a case distinction.

In the case $x>y$ we get an $\varepsilon$-interleaving if $\varepsilon > x-y$, but not if $\varepsilon < x-y$, by symmetry.
Combining both cases we get 
\begin{align*}
d_I(M,N)= |x-y| \ . & \qedhere
\end{align*}
\end{proof}
\end{lem}

As extended pseudometric $d_I$ induces a topology $\Tid$ on $\Sp \RRR$, in contrast to the topology $\Tcl$ coming from the closure operation, with basis the open balls
\begin{align*}
B_{\varepsilon}^{d_I}(k_I)= \left\lbrace L \in \idl \RRR | d_I(k_I,k_L) < \varepsilon \right\rbrace =  \Delta \left(\principal {(x-\varepsilon)}, \qrincipal {(x + \varepsilon)} \right)
\end{align*}
for every $\varepsilon \geq 0$ and $I = \principal x$ or $I = \qrincipal x \in \idl \RRR$.
We observe that the interleaving distance cannot distinguish between modules of the form $k_{\qrincipal x}$ and $k_{\principal x}$ in the spectrum $\Sp \RRR$, and therefore the induced topological space is not Hausdorff.
Note that $d_I(k_{\RRR}, M) = \infty$ for all $M \in \Sp \RRR \setminus \lbrace k_{\RRR} \rbrace$, so $B_{\varepsilon}^{d_I}(k_{\RRR})= \lbrace k_{\RRR} \rbrace$.

The open ball $B_{\varepsilon}^{d_I}(k_I)$ is the complement of the closed set 
\begin{align*}
\Delta \left( \left( -\infty, \principal{(x-\varepsilon)} \right] \right) \cup  \Delta \left( \left[ \principal{(x+\varepsilon)}, \infty \right] \right)
\end{align*}
if $I \neq \RRR$, so these open balls in the interleaving topology are also Ziegler open.
But for any $\varepsilon > 0$ the complement of $B^{d_I}_{\varepsilon}(k_{\RRR})$ is $\Delta \left( (-\infty, \RRR) \right)$, so  $\lbrace k_{\RRR} \rbrace$ is not Ziegler closed. 
Hence, the topology induced by the interleaving distance is only refined by the Ziegler topology if we exclude the point at infinity $k_{\RRR}$.
\begin{cor}
The restricted identity morphism
\begin{gather*}
(\Sp \RRR \setminus  \lbrace k_{\RRR} \rbrace, \Tcl | _{\Sp \RRR \setminus  \lbrace k_{\RRR} \rbrace}) \rightarrow (\Sp \RRR, \Tid)
\end{gather*}
is continuous.
\end{cor}

\bibliography{Literature}
\bibliographystyle{amsalpha}
\end{document}